\documentclass[12pt,letterpaper]{amsart}
\usepackage[utf8x]{inputenc}

\usepackage{amsfonts,amssymb,amstext,amsmath,amscd,latexsym,amsthm}

\usepackage{answers}
\Newassociation{sol}{Solution}{ans}

\usepackage{pstricks}
\usepackage{pst-plot}

\theoremstyle{plain}
\newtheorem{theorem}{Theorem}
\theoremstyle{definition}
\newtheorem{definition}{Definition}

\theoremstyle{remark}
\newtheorem{example}{Example}
\newtheorem{exercise}{Exercise}
\newtheorem{remark}{Remark}

\date{}
\author{J. R. Arteaga, M. Malakhaltsev}
\title[$G$-structures and differential equations]{Ideas of E.~Cartan and S.~Lie in modern geometry: $G$-structures and differential equations. 
Lecture 3}
\begin{document}
\Opensolutionfile{ans}[Lecture3_answers]
\maketitle

\fbox{\fbox{\parbox{5.5in}{
\textbf{Problem:}\\
How to make Cartan reduction in a particular case.
}}}

\vspace{1cm}

We will show how to find a differential invariant of a $G$-structure by example of a contact $2$-distribution in $\mathbb{R}^{3}$, using the Cartan reduction method step by step.

\section*{Two-dimensional distribution $\Delta$ in $M^{3}$} 
\label{sec:two-dimensional_distribution_in_R3}

\begin{definition}
Let $M$ a three-dimensional manifold. A $2$-dimensional distribution $ \Delta $ on $M$ is an assignment of  a plane to each point of  $M$, i.e.  $\Delta$ is a sub-bundle of the tangent bundle $TM$. 
The assignment is smooth in sense that in a neighborhood $U$ of each point $p \in M$ there are vector fields $\{X_{1}, X_{2}\}$ such that:
\begin{enumerate}
\item 
For any point  $q \in U$ the vectors $X_{1}(q)$, $X_{2}(q)$ are linear independent.
\item
For each point $ q\in U $,  $\Delta_{q}$ is the plane  spanned by the two vectors $X_{1}(q)$, $X_{2}(q)$.
\end{enumerate}
\end{definition}

\subsection*{Sub-riemannian surface $\mathcal{S}$ on $M$}
 
Let $M$ a three-dimensional manifold. 
A distribution $\Delta$ is called \textit{integrable} or \textit{holonomic distribution}  if for each point $ p\in M$ there exists  a surface $\Sigma$ passing through $p$ which is tangent to $\Delta$: $T_q\Sigma = \Delta(q)$ for each $q \in \Sigma$. 
A $ 2 $-distribution $ \Delta $ on $M$ is holonomic if the \emph{commutator }of vector fields $ X_{1} $ and $ X_{2} $, 
\begin{equation}
[X_{1}, X_{2}]^{i} = 
X_{1}^{s} \dfrac{\partial X_{2}^{i}}{\partial x^{s}} -
X_{2}^{s} \dfrac{\partial X_{1}^{i}}{\partial x^{s}}
\end{equation}
generating $\Delta $ belong to $ \Delta $, i.e. $ [X_{1}, X_{2}](p) \in \Delta(p)$. In this case the family of these surfaces form  a \emph{foliation} of $M$. 

In another case, when $  [X_{1}, X_{2}] \notin \Delta  $ for all points in $M$, we say the distribution is \textit{non-integrable} or \textit{non-holonomic}. 

\begin{definition}
A sub-riemannian surface in $M$ is  a non-holonomic distribution $\Delta$ with a scalar product $\langle \cdot , \cdot \rangle$ on $\Delta(p)$ for each $p \in M$. We denote by $\mathcal{S} = \left(M, \left(\Delta, \langle \cdot , \cdot \rangle \right) \right)$ a sub-riemannian surface in $M$.
\end{definition}

In this lecture we will show how to use the Cartan ideas for to  find some invariants of a sub-riemannian surface $\mathcal{S}$. 

\begin{example}[Heisenberg distribution]\label{Exa:example1}
Let  $\mathcal{S} = \left(\mathbb{R}^{3}, \Delta, \langle \cdot , \cdot \rangle \right)$ be the sub-riemannian surface where
\begin{equation}
\Delta = span \{X_{1}, X_{2}\}, \quad 
\begin{cases}
X_{1}=(1,0,-y) = \dfrac{\partial}{\partial x} - y \dfrac{\partial}{\partial z}\\
X_{2}=(0,1,x) = \dfrac{\partial}{\partial y} + x \dfrac{\partial}{\partial z}\\
\end{cases}
\end{equation}
and the scalar product is the induced from $ \mathbb{R}^{3} $. This distribution is non-holonomic because
\begin{equation}
\begin{split}
[X_{1}, X_{2}] & = Y \\
Y & = X_{3}^{s}\dfrac{\partial}{\partial x^{s}}\\
\text{where}\quad   Y^{s} & = X_{1}^{k}\dfrac{\partial X_{2}^{s}}{\partial x^{k}} -
X_{2}^{k}\dfrac{\partial X_{1}^{s}}{\partial x^{k}}\\
\therefore \quad Y &= (0,0,2) = 2 \dfrac{\partial}{\partial z} \notin \Delta
\end{split}
\end{equation}
\end{example}

\begin{exercise}\label{Exe:exercise1}
Let  $\mathcal{S} = \left(\mathbb{R}^{3}, \Delta, \langle \cdot , \cdot \rangle \right)$ be the sub-riemannian surface where
\begin{description}
\item[(a)]
$ X_{1}=(1, 0, x_{2})  $, $ X_{2}=(0, 1, x_{1})  $ in $ \mathbb{R}^{3} $.
\item[(b)]
$ X_{1}=(1, 0, -x_{2})  $, $ X_{2}=(0, 1, 0) $ in $ \mathbb{R}^{3} $.
\end{description}
and the scalar product is the induced from $ \mathbb{R}^{3} $.
Is the distribution $ \Delta = span \{X_{1}, X_{2}\} $ holonomic or not?
\begin{sol}
\begin{description}
\item[(a)]
Holonomic distribution because $[X_{1}, X_{2}]=0$.
\item[(b)]
Non-holonomic distribution because $[X_{1}, X_{2}]=(0,0,1)$. This distribution is called Cartan distribution.
\end{description}
\end{sol}
\end{exercise}

\section*{Cartan reduction for a sub-Riemannian surfaces $\mathcal{S}$ in $\mathbb{R}^{3}$}\label{sec:Cartan_reduction}

Let $M$ be a three-dimensional manifold. Consider a sub-riemannian surface $\mathcal{S} = \left(M, \Delta, \langle \cdot, \cdot \rangle \right)$. 
For each $p\in M$ we always can take a local orthonormal frame field $(e_{1}, e_{2})$ of $\Delta$. If we take $e_{3}=[e_{1}, e_{2}]$ then $(e_{1}, e_{2}, e_{3})$ is a local frame field of $M$.

Let $B(M)$ be the bundle of positively oriented coframes of $M$.

\subsection*{\underline{First step:} Adapting a coframe}
\begin{definition}
Given an oriented sub-riemannian surface $\mathcal{S} = (\Delta, \langle \cdot , \cdot \rangle)$ on a three-dimensional manifold $M$,
we say that a co-frame $\eta = \left(\eta^{1}, \eta^{2}, \eta^{3} \right)$ of $B(M)$, is \emph{adapted} to $\Delta$ if for all $ p\in M$,
\begin{itemize}
\item[(1)]%
$(\eta^1 |_{\Delta_p}, \eta^2 |_{\Delta_p})$ is a positively oriented co-frame of $\Delta_p$;
\item[(2)]%
$\eta^3 (W) = 0$ for any $W \in \Delta_p$; 
\item[(3)]%
$ \langle W , W \rangle = [\eta^{1}(W)]^{2} + [\eta^{2}(W)]^2$ for any $W \in \Delta_p$.
\end{itemize}
\end{definition}

\begin{definition}[Sub-bundle $B_{0}$]
The sub-bundle $B_0 \subset B$ consists  of all co-frames adapted to $\mathcal{S}$. 
The subgroup of matrices of $GL(3)$ which transform adapted coframes into adapted coframes 
is
\begin{equation}
G_0 = \left \{
\left.
\left(
\begin{array}{ccc}
\cos\varphi_{1} & -\sin\varphi_{1} & \varphi_{2}\\
\sin\varphi_{1} & \cos\varphi_{1} & \varphi_{3}\\
0 & 0 & \varphi_{4}  
\end{array}
\right)
\ \right| \varphi_{4} \ne 0
\right\}
\label{eq:G_0}
\end{equation}
\end{definition}
With this construction we show that any sub-riemannian surface  $(\Delta,\langle \cdot,\cdot \rangle)$ defines a $G_0$-structure on $M$.

\begin{remark}
The quantities $\varphi_{i}$ ($i\in \{1,2,3,4\}$) are real variables and can be considered as coordinates on $G_0$, so the dimension of $G_0$ is $4$, and then the dimension of $B_0$ is $7$. 

Since $G_{0}$ is a Lie group with identity element $I$,  one can construct the associated Lie algebra $\mathfrak{g}_{0}$ as the tangent space of $G_{0}$ in $I$. 
As coordinates of $I$ are $\varphi_1 = \varphi_2 = \varphi_3 = 0$ and $\varphi_4=1$, taking the tangent vectors to the coordinate curves, we obtain that 
the Lie algebra $\mathfrak{g}_{0}$ associated to the Lie group $G_{0}$ is the set of matrices of type
\begin{equation}
\left (
\begin{array}{ccc}
0 & \alpha_{1} & \alpha_{2} \\
-\alpha_{1} & 0 & \alpha_{3} \\
0 & 0 & \alpha_{4}
\end{array}
\right )
\end{equation}
\end{remark}

\begin{example}\label{Exa:example2}
Find an adapted frame and coframe for  the Heisenberg distribution of example \eqref{Exa:example1},
\begin{equation}
X_{1}=(1, 0, -y), \quad  X_{2}=(0, 1, x)
\end{equation}
This distribution is defined in whole $ \mathbb{R}^{3} $ and is given by the vector fields $ \{ X_{1}, X_{2}\} $. Using Gram-Schmidt algorithm we found one orthonormal frame for $ \Delta $, $ (e_{1}, e_{2}) $ and complete it with $ e_{3} = [e_{1}, e_{2}]$ to obtain the frame $e=(e_{1}, e_{2}, e_{3})  $ for $ \mathbb{R}^{3} $. So $e=(e_{1}, e_{2}, e_{3})$ is an adapted frame for $\mathbb{R}^{3}$,
\begin{equation}\nonumber
\begin{cases}
e_{1}=
\dfrac{1}{\sqrt{1+y^{2}}}\dfrac{\partial }{\partial x}-
\dfrac{y}{\sqrt{1+y^{2}}}\dfrac{\partial }{\partial z}\\
e_{2}=
\dfrac{xy}{\sqrt{1+y^{2}}\sqrt{1+x^{2}+y^{2}}}\dfrac{\partial }{\partial x}+
\dfrac{\sqrt{1+y^{2}}}{\sqrt{1+x^{2}+y^{2}}}\dfrac{\partial }{\partial y}+
\dfrac{x}{\sqrt{1+y^{2}}\sqrt{1+x^{2}+y^{2}}}\dfrac{\partial }{\partial z}\\
e_{3}=
\dfrac{y}{\left (1+x^{2}+y^{2}\right )^{3/2}}\dfrac{\partial }{\partial x}- 
\dfrac{x}{\left (1+x^{2}+y^{2}\right )^{3/2}}\dfrac{\partial }{\partial y}-
\dfrac{2+3x^{2}+3y^{2}}{\left (1+x^{2}+y^{2}\right )^{3/2}}\dfrac{\partial }{\partial z}
\end{cases}
\end{equation}

The dual co-frame is
\begin{equation}\nonumber
\begin{cases}
\eta^{1}=\dfrac{\left (2+3y^{2}\right )dx}{2\sqrt{1+y^{2}}}-\dfrac{3xydy}{2\sqrt{1+y^{2}}} +\dfrac{ydz}{2\sqrt{1+y^{2}}}\\
\eta^{2}=-\dfrac{xydx}{2\sqrt{1+y^{2}}\sqrt{1+x^{2}+y^{2}}}+\dfrac{\left (2+3x^{2}+2y^{2}\right )dy}{2\sqrt{1+y^{2}}\sqrt{1+x^{2}+y^{2}}} -\dfrac{xdz}{2\sqrt{1+y^{2}}\sqrt{1+x^{2}+y^{2}}}\\
\eta^{3}=-\dfrac{y}{2}\sqrt{1+x^{2}+y^{2}} dx+\dfrac{x}{2}\sqrt{1+x^{2}+y^{2}} dy-\dfrac{1}{2}\sqrt{1+x^{2}+y^{2}} dz
\end{cases}
\end{equation}
that is an adapted co-frame to $\Delta$.
\end{example}

\begin{remark}
We can write that the Heisenberg distribution is $\Delta = \ker (\eta^{3})$.
\end{remark}

\begin{exercise}\label{Exe:exercise2}
Find an adapted frame and coframe for  the Cartan distribution of exercise \eqref{Exe:exercise1}
\begin{equation}
X_{1}=(1, 0, -y), \quad  X_{2}=(0, 1, 0) 
\end{equation}

\begin{sol}
Adapted frame
\begin{equation}
e = 
\left(
\begin{matrix}
\dfrac{1}{\sqrt{1+y^{2}}} & 0 & \dfrac{y}{\sqrt{1+y^{2}}}\\
0 & 1 & 0\\
-\dfrac{y}{\sqrt{1+y^{2}}}  & 0 & \dfrac{1}{\sqrt{1+y^{2}}} 
\end{matrix}
\right)
\end{equation}
The adapted co-frame $ \eta = (\eta^{1}, \eta^{2}, \eta^{3})$ for $ \Delta $ is:
\begin{equation}
\begin{cases}
\eta^{1} = \dfrac{1}{\sqrt{1+y^{2}}} dx - \dfrac{y}{\sqrt{1+y^{2}}}dz \\ 
\eta^{2} =dy\\
\eta^{3} =  -\dfrac{y}{\sqrt{1+y^{2}}} dx + \dfrac{1}{\sqrt{1+y^{2}}}dz 
\end{cases}
\end{equation}
\end{sol}
\end{exercise}

\begin{definition}[Tautological forms]
The $1$-forms $\theta^{i}$ on $B_0$ such that $\theta^i(\eta^a)(X) = \eta^i(d\pi(X))$ are called tautological forms.

An adapted co-frame $\eta = (\eta^{1}, \eta^{2}, \eta^{3})$ defined on a neighborhood $U \subset M$
defines a trivialization 
\begin{equation}
U \times G_0 \leftrightarrow B_0|_U , \quad (x,g) \leftrightarrow g^{-1}\eta_x.
\end{equation}
In terms of this trivialization, the tautological forms can be written as 
\begin{equation}
\theta_{(x,g)} = g^{-1} (\eta_x \circ d\pi) 
\end{equation}
\end{definition}

\subsection*{Derivation equations}
\begin{theorem}
Let $\mathcal{S}=(\Delta , \langle \cdot , \cdot \rangle)$ be a sub-riemannian surface in $M$. 
Let $B(M)$ be the principal bundle of positive oriented co-frames on $M$ and $B_{0}(M)$ the principal sub-bundle of $B(M)$ with the group $G_{0}$ defined in \eqref{eq:G_0} consisting of co-frames adapted to $\mathcal{S}$. 
The exterior derivatives of the tautological forms can be written as follows:
\begin{equation}
\left (
\begin{array}{c}
d\theta^{1}\\
d\theta^{2}\\
d\theta^{3}
\end{array}
\right )
=
\left (
\begin{array}{ccc}
0 & \alpha_{1} & \alpha_{2} \\
-\alpha_{1} & 0 & \alpha_{3} \\
0 & 0 & \alpha_{4} 
\end{array}
\right )
\wedge
\left (
\begin{array}{c}
\theta^{1}\\
\theta^{2}\\
\theta^{3}
\end{array}
\right )
+
\left (
\begin{array}{ccc}
T_{23}^{1} & T_{31}^{1} & T_{12}^{1} \\
T_{23}^{2} & T_{31}^{2} & T_{12}^{2} \\
T_{23}^{3} & T_{31}^{3} & T_{12}^{3} 
\end{array}
\right )
\left (
\begin{array}{c}
\theta^{2}\wedge\theta^{3}\\
\theta^{3}\wedge\theta^{1}\\
\theta^{1}\wedge\theta^{2}
\end{array}
\right )
\label{Th:theorema_1}
\end{equation}
or in contracted form,
\begin{equation}
\boxed{
d \theta ^{i} = \omega^{i}_{s} \wedge \theta^{s} + T^{i}_{ab} \theta^{a}\wedge\theta^{b}
}
\end{equation}
\end{theorem}

\begin{remark}
\begin{enumerate}
\item
These equations \eqref{Th:theorema_1} are a part of all structure equations. They represent only the equations containing the derivatives of the tautological forms $\theta^{i}$ in terms of themselves. Recall that the cotangent space $T^{*}B_{0}$ has dimension $7$. One co-frame of $T^{*}B_{0}$ is:
\begin{equation}
\{ \theta^{1}, \theta^{2}, \theta^{3}, d\alpha_{1}, d\alpha_{2}, d\alpha_{3}, d\alpha_{4}\}
\end{equation}
\item
The $1$-form $\omega = \omega_{a}^{b}$ is a  pseudoconnection form. 

\item
The coefficients $T^{i}_{ab}$ are called the torsion coefficients.
\item
We would like to construct the invariants from the components of the connection form $\omega$ and the torsion $T$. 

Obviously these components in general case are functions, they depends of the points $(x,g)$ on a fibre $T_{x}B$. If all of them were constant we would finish the problem, however in the general case they are not. For this reason we have to continue.
\end{enumerate}
\end{remark}

\subsubsection*{Contact distribution}

A distribution $\Delta$ on $M$ is a contact distribution if for all $p\in M$ the plane $\Delta(p)$ is given by the zeros of a $1$-form $\eta^{3}$. Is clear that they will also given by the zeros of $\lambda \eta^{3}$. Thus, $\{\lambda \eta^{3}\}$ all give same the same $\Delta(p)$. The Heisenberg distribution and the Cartan distribution are both contact distribution.

The property that a contact distribution $\Delta$ in $\mathbb{R}^{3}$  is non-integrable if
\begin{equation}
d \eta^{3} \wedge \eta^{3} \ne 0
\end{equation}

This property also can be formulated in terms of the tautological form: a contact distribution is a non-integrable if
\begin{equation}
d \theta^{3} \wedge \theta^{3} \neq 0
\end{equation}
In our case, using the equation \eqref{Th:theorema_1}, for a contact distribution in $\mathbb{R}^{3}$ we have,
\begin{equation}
d \theta^{3} \wedge \theta^{3} = T_{12}^{3}\theta^{1}\wedge\theta^{2}\wedge\theta^{3}
\end{equation}

\begin{example}\label{Exa:example_3}
For the Heisenberg distribution treated in the example \ref{Exa:example2} we have,
\begin{equation}
d \theta^{3} \wedge \theta^{3} = -2 \theta^{1}\wedge\theta^{2}\wedge\theta^{3}
\end{equation}
That is, $T^{3}_{12}=-2$
\end{example}

\begin{exercise}\label{Exe:exercise_3}
For the Cartan distribution find $T^{3}_{12}$.
\begin{sol}
$T^{3}_{12}=-1$
\end{sol}
\end{exercise}

\subsection*{How the component $T^{3}_{12}$ change under the action of $G_{0}$?}
Let us denote the action of $G_{0}$ on $B^{0}$ by $R_{g}$ where $g\in G_{0}$. Thus the components of  $\theta^{i}$ under action of $G_{0}$ change by the following rule,
\begin{equation}
R_{g}\theta = g^{-1}\theta
\end{equation}
This imply that
$R_{g}\theta^{3} = g^{-1}\theta^{3}$, where $g\in G_{0}$ defined in equation \eqref{eq:G_0}.  Therefore $R^{*}_{g}d\theta^{3} = (\varphi_{4})^{-1}d\theta^{3}$, and this imply that
\begin{equation}
R^{*}_{g}T^{3}_{12} = (\varphi_{4})^{-1}T^{3}_{12} 
\label{eq:3_1}
\end{equation}
For a contact distribution $T^{3}_{12} \ne 0$, and from \eqref{eq:3_1} it follows that  we can take a subbundle $B_1 \subset B_0$ 
with the property that $T^3_{12}=1$. The structure group $G_1$ of $B_1$ is   
\begin{equation}
\boxed{
G_1 = \left \{
\left(
\begin{array}{ccc}
\cos\varphi_{1} & -\sin\varphi_{1} & \varphi_{2}\\
\sin\varphi_{1} & \cos\varphi_{1} & \varphi_{3}\\
0 & 0 & 1
\end{array}
\right)
\right\}
\label{eq:G_0_new}
}
\end{equation}
and the Lie algebra $\mathfrak{g}_{0}$ associated to $G_{0}$ is the algebra of matrices,
\begin{equation}
\boxed{
\mathfrak{g}_{1} = \left \{
\left (
\begin{array}{ccc}
0 & \alpha_{1} & \alpha_{2} \\
-\alpha_{1} & 0 & \alpha_{3} \\
0 & 0 & 0
\end{array}
\right )
\right\}
}
\end{equation}
Thus we have reduced the structure group $G_{0}$ to $G_1$, and the dimension of $B_{1}$ is $6$, because the dimension of the vertical space, isomorphic to $G_1$, is $3$.

\subsection*{\underline{Second step}: Reducing the structure group $G_{0}$ as much as possible}

The idea is continue reducing the structure group $G_{0}$. Now we have a principal  sub-bundle $(B_{1}, G_{1})$ of $(B_{0}, G_{0})$ defined by   
\begin{eqnarray}
&& B_{1}=\left \{ \eta = \left (\eta^{1}, \eta^{2}, \eta^{3} \right ) \in B_{0}\mid T^{3}_{12}(\eta)=1 \right \} \nonumber\\
&&
G_{1}=\left \{
\left (
\begin{array}{ccc}
\cos\varphi_{1} & -\sin\varphi_{1} & \varphi_{2}\\
\sin\varphi_{1} & \cos\varphi_{1} & \varphi_{3} \\
0 & 0 & 1  
\end{array}
\right )=
\left (
\begin{array}{cc}
A & B \\
0 & 1
\end{array}
\right )
\right \}
\nonumber
\end{eqnarray}
The structure equations are,
\begin{equation}
\left (
\begin{array}{c}
d\theta^{1}\\
d\theta^{2}\\
d\theta^{3}
\end{array}
\right )
=
\left (
\begin{array}{ccc}
0 & \alpha &\beta \\
-\alpha & 0 & \gamma \\
0 & 0 & 0 
\end{array}
\right )
\wedge
\left (
\begin{array}{c}
\theta^{1}\\
\theta^{2}\\
\theta^{3}
\end{array}
\right )
+
\left (
\begin{array}{ccc}
T_{23}^{1} & T_{31}^{1} &T_{12}^{1} \\
T_{23}^{2} & T_{31}^{2} &T_{12}^{2} \\
T_{23}^{3} & T_{31}^{3} & 1 
\end{array}
\right )
\left (
\begin{array}{c}
\theta^{2}\wedge\theta^{3}\\
\theta^{3}\wedge\theta^{1}\\
\theta^{1}\wedge\theta^{2}
\end{array}
\right )\nonumber
\end{equation}

\subsection*{How the components $T^{3}_{23}$ and $T^{3}_{31}$ are changed under the action of $G_{1}$?}
Since the right action $R_{g}$ is defined by $R_{g} \theta = g^{-1}\theta$ where $g\in G_{1}$, then the components $T^{3}_{23}$ and $T^{3}_{31}$ under action of $G_{1}$ change by the following rule:
\begin{equation}\nonumber
\boxed{
R_{g}^{*}
\left (
\begin{array}{c}
T^{3}_{23}\\
T^{3}_{31}
\end{array}
\right )
=
A^{-1}
\left (
\begin{array}{c}
T^{3}_{23} -\varphi_{2}\\
T^{3}_{31}-\varphi_{3}
\end{array}
\right )
}
\end{equation}

So, we can define another principal sub-bundle $(B_{2}, G_{2})$ as follows, 

\begin{eqnarray}
&& B_{2}=\left \{ \eta = \left (\eta^{1}, \eta^{2}, \eta^{3} \right ) \in B_{1}\mid T^{3}_{23}(\eta)=T^{3}_{31}(\eta)=0 \right \} \nonumber\\
&&
G_{2}=\left \{
\left (
\begin{array}{ccc}
\cos\varphi & -\sin\varphi & 0\\
\sin\varphi & \cos\varphi & 0 \\
0 & 0 & 1  
\end{array}
\right )=
\left (
\begin{array}{cc}
A & 0 \\
0 & 1
\end{array}
\right )
\right \}
\nonumber
\end{eqnarray}
and the Lie algebra $\mathfrak{g}_{2}$ associated to $G_{2}$ is

\begin{equation}
\mathfrak{g}_{2} = \left \{
\left (
\begin{array}{ccc}
0 & \alpha_{1} & 0\\
-\alpha_{1} & 0 & 0 \\
0 & 0 & 0
\end{array}
\right )
\right\}
\end{equation}

Then the structure equations are,
\begin{equation}
\left (
\begin{array}{c}
d\theta^{1}\\
d\theta^{2}\\
d\theta^{3}
\end{array}
\right )
=
\left (
\begin{array}{ccc}
0 & \alpha & 0 \\
-\alpha & 0 & 0 \\
0 & 0 & 0 
\end{array}
\right )
\wedge
\left (
\begin{array}{c}
\theta^{1}\\
\theta^{2}\\
\theta^{3}
\end{array}
\right )
+
\left (
\begin{array}{ccc}
T_{23}^{1} & T_{31}^{1} &T_{12}^{1} \\
T_{23}^{2} & T_{31}^{2} &T_{12}^{2} \\
0 & 0 & 1 
\end{array}
\right )
\left (
\begin{array}{c}
\theta^{2}\wedge\theta^{3}\\
\theta^{3}\wedge\theta^{1}\\
\theta^{1}\wedge\theta^{2}
\end{array}
\right )
\label{Eq:structure_equation_B2}
\end{equation}


\subsection*{\underline{Third step}: Finding invariants}

We have reduced $G_{0}$ to a minimum subgroup $G_{2}$ imposing more and more conditions for the adapted co-frame $\eta$. In the ideal case  we must find the unique frame of the $G$-structure in order to find finally invariants of $G$-structurs.

\subsubsection*{How the components of $\omega$ and $T$ in the structure equations are changed?} 
If we take another connection form $\omega'$, then the torsion map also changes and we have $T'$ such that
\begin{equation}
d\theta^{i} = \omega'^{i}_{m} \wedge \theta^{m} + T'^{i}_{lm} \theta^{l} \wedge \theta^{m} \nonumber \\
\end{equation}


But $\omega$ and $\omega'$ are both connections on $B_{2}$, i.e.  both are  elements  in the space of all smooth $1$-forms on $B_{2}$ with values in the Lie algebra $\mathfrak{g}_{2}$, $ \Lambda^{1}(B_{2}, \mathfrak{g}_{2})$,

If $\sigma$ is a fundamental vector field and $\mu$ is the difference $\mu = \omega' - \omega$ we have
\begin{equation}
\omega (\sigma (a))=a, \omega' (\sigma (a))=a \quad \text{where}\quad a\in \mathfrak{g}_{2}
\end{equation}
and
\begin{equation}
\mu (\sigma(a))=0
\end{equation}
Therefore $\mu$ does vanish on the vertical subbundle $V$. 
So we have,
\begin{equation}
\boxed{
\mu = \mu^{i}_{js}\theta^{s} \Rightarrow
\omega'^{i}_{m} = \omega^{i}_{m} + \mu^{i}_{ms}\theta^{s}
}
\end{equation}
Therefore,
\begin{equation}
\begin{split}
d\theta^{i} &= \omega'^{i}_{m} \wedge \theta^{m} + T'^{i}_{lm} \theta^{l} \wedge \theta^{m} \nonumber \\
&= \left (\omega^{i}_{m} + \mu^{i}_{ml}\theta^{l}\right )\wedge \theta^{m} + T'^{i}_{lm} \theta^{l} \wedge \theta^{m} \nonumber \\
\end{split}
\end{equation}
and
\begin{equation}
\begin{split}
d\theta^{i} &= \omega'^{i}_{m} \wedge \theta^{m} + T'^{i}_{lm} \theta^{l} \wedge \theta^{m} \\
&= \left (\omega^{i}_{m} + \mu^{i}_{ml}\theta^{l}\right )\wedge \theta^{m} + T'^{i}_{lm} \theta^{l} \wedge \theta^{m}  \\
&= \omega^{i}_{s}\wedge \theta^{s} + \left (T'^{i}_{lm} - \mu^{i}_{[lm]}\right ) \theta^{l} \wedge \theta^{m} =
\omega^{i}_{s}\wedge \theta^{s} + T^{i}_{lm} \theta^{l} \wedge \theta^{m}\\
&
\boxed{
T'^{i}_{lm} = T^{i}_{lm} + \mu^{i}_{[lm]}
}
, \quad 
\text{where}
\quad
\mu^{i}_{[lm]}=
A(\mu^{i}_{lm})=\mu^{i}_{ml}
\end{split}
\end{equation}
where $A$ is the operator of alternation with restpect to the lower indices. 

\subsubsection*{Finding invariants}
If we take 
\begin{equation}
\alpha' = \alpha + T^{1}_{12}\theta^{1}+T^{2}_{12}\theta^{2}-\dfrac{1}{2}\left (T^{2}_{31}+T^{1}_{23}\right )\theta^{3}
\end{equation}
and replace it in \eqref{Eq:structure_equation_B2} we obtain,
\begin{equation}
T'^{1}_{12}=T'^{2}_{12}=0,\quad \text{and}
\quad
T'^{1}_{23}=-T'^{2}_{31}
\end{equation}
and differentiating $d(\theta^{3})$ we obtain,
\begin{equation}
d(d\theta^{3}) =0\Longrightarrow 
T'^{1}_{31}=T'^{2}_{32}
\end{equation}

If we denote $a_{1}=T^{1}_{23}, a_{2}=T^{1}_{31} $, then the  structure equation will be reduce to,

\begin{equation}
\begin{split}
&
\left (
\begin{array}{c}
d\theta^{1}\\
d\theta^{2}\\
d\theta^{3}
\end{array}
\right )
=
\left (
\begin{array}{ccc}
0 & \alpha & 0 \\
-\alpha & 0 & 0 \\
0 & 0 & 0 
\end{array}
\right )
\wedge
\left (
\begin{array}{c}
\theta^{1}\\
\theta^{2}\\
\theta^{3}
\end{array}
\right )
+
\left (
\begin{array}{ccc}
a_{1} & a_{2} & 0 \\
a_{2} & -a_{1} & 0 \\
0 & 0 & 1 
\end{array}
\right )
\left (
\begin{array}{c}
\theta^{2}\wedge\theta^{3}\\
\theta^{3}\wedge\theta^{1}\\
\theta^{1}\wedge\theta^{2}
\end{array}
\right )
\end{split}
\end{equation}
One differential invariant of a contact distribution on $\mathbb{R}^{3}$ is,
\begin{equation}
\mathcal{M} = (a_{1})^{2} + (a_{2})^{2}
\end{equation} 

\begin{example}
Let us consider the Heisenberg distribution defined on $\mathbb{R}^{3}$ treated in the examples \eqref{Exa:example2}, and \eqref{Exa:example_3}, 
\begin{equation}
\eta^{3} = ydx -xdy +dz
\end{equation}
For this contact distribution, using the Cartan reduction method, we have that between another differential invariants has the invariant 
\begin{equation}
\mathcal{M} = \dfrac{9}{4}\dfrac{(x^{2}+y^{2})^{2}}{(1+x^{2}+y^{2})^{4}}
\end{equation} 
These calculations we obtained using Maple.
\end{example}

\begin{exercise}
Calculated the invariant $\mathcal{M}$ for the Cartan distribution treated in the second and third exercise. Is the Cartan distribution equivalent to the Heisenberg distribution?.
\begin{sol}
\begin{equation}
\mathcal{M} = \dfrac{1}{4}\dfrac{(2y^{2}-1)^{2}}{(1+y^{2})^{4}}
\end{equation}
The Cartan and Heisenberg distributions are not equivalents. 
\end{sol} 
\end{exercise}

\section*{Summary of Lecture 3}
\begin{enumerate}
\item
The Cartan reduction method is a tool in the modern Differential Geometry in order to determine if two geometrical structures are equivalent up to a diffeomorphism. 
We demonstrated this method by an example of a contact $2$-distribution with a metric $\mathcal{S} = (\Delta, \langle \cdot , \cdot \rangle)$ in a three-dimensional manifold $M$ and found a differential invariant for this geometrical structure.
\item
The method consist of three steps: 
\begin{enumerate}
\item
At the first step one defines an adapted coframe for the distribution $\mathcal{S}$ and constructs a sub-bundle $(B_{0}, G_{0})$ of the principal $GL(3)^{+}$-bundle of all oriented positively coframes on $M$. 
We write the structure equations which express the exterior derivatives of tautological forms in terms of themselves and a connection form.  
\item
At the second step we reduce the group $G_{0}$ as much as possible adding new conditions for the adapted coframe such that the structure group becomes smaller and smaller. Thus we get a sub-bundle  $(B_{2}, G_{2})$.
\item
At the third step we construct invariants from the torsion coefficients. 

\end{enumerate}

\end{enumerate}

\Closesolutionfile{ans}
\section*{Answers to exercises}
\input{Lecture3_answers}


\end{document}